\documentclass[journal]{IEEEtran}

\usepackage{amsmath,amssymb,bm}
\usepackage{graphicx}
\usepackage{booktabs}
\usepackage{array}
\usepackage{multirow}
\usepackage{cite}
\usepackage{algorithm}
\usepackage{algpseudocode}
\usepackage[caption=false,font=footnotesize]{subfig}
\usepackage{url}
\usepackage{xcolor}
\usepackage[hidelinks]{hyperref}

\graphicspath{{figures/}}
\newcommand{\R}{\mathbb{R}}

\newcommand{\SOC}{\mathrm{SOC}}
\newcommand{\argmin}{\mathop{\mathrm{arg\,min}}}
\newcommand{\calI}{\mathcal{I}}
\newcommand{\calS}{\mathcal{S}}
\newcommand{\calB}{\mathcal{B}}

\newcommand{\revred}[1]{#1}

\begin{document}

\title{Battery-Swapping Station Operation Under Forecast Uncertainty: A Scenario-Based Stochastic MPC Framework}

\author{\begin{tabular}{c}
Zhiyuan Guo\textsuperscript{1}, Siyang Gao\textsuperscript{1}, Zhichao Chen\textsuperscript{1}, Zhankun Sun\textsuperscript{2} 
, and Jiaze Ma*\textsuperscript{1}\\[0.4ex] \small\textsuperscript{1}Department of Systems Engineering, City University of Hong Kong, Kowloon, Hong Kong SAR, China\\
\small\textsuperscript{2}Department of Decision Analytics and Operations, City University of Hong Kong, Kowloon, Hong Kong SAR, China\\
\thanks{Corresponding author: Jiaze Ma (e-mail: jiazema@cityu.edu.hk).}
\end{tabular}}

\maketitle

\begin{abstract}
Battery-swapping stations (BSSs) can shorten electric-vehicle energy replenishment while using centrally managed battery inventories as flexible grid-connected storage. Realizing both benefits requires the station to schedule charging, grid discharge, and swapping service before future customer demand and electricity prices are known. This paper develops a forecast-aware rolling-horizon operating framework for this problem. A lightweight DLinear model predicts 24-hour price and demand trajectories, Stein variational gradient descent quantifies their uncertainty through representative scenarios, and a two-stage stochastic model predictive controller converts those scenarios into station decisions. The controller accounts for service shortfall, terminal readiness, a protected service buffer, and electrochemical degradation without assuming perfect future information. The application contribution is an implementable controller that coordinates the station's mobility-service and energy-storage roles. The methodological contribution is a modular forecast-to-control interface that separates the operational value of mean-forecast accuracy from that of uncertainty representation. In a 120-day closed-loop evaluation, DLinear-SVGD SMPC achieves the lowest cost among the implementable controllers. Relative to deterministic DLinear MPC, it reduces final cost by 1.2\% and service-shortfall hours by 80.7\%, with 99.10\% of the evaluated hours free of shortfall.
\end{abstract}

{\renewcommand{\abstractname}{Note to Practitioners}
\begin{abstract}
A battery-swapping station keeps a pool of charged batteries so that an arriving driver can exchange a depleted pack instead of waiting for the vehicle to recharge. The same battery pool can also behave as station-based energy storage: the operator can buy electricity when it is inexpensive, postpone charging during high-price periods, and sell energy when doing so will not compromise customer service. The practical difficulty is that future swap arrivals and electricity prices are uncertain, yet charging and energy-sale decisions must be made in advance. This paper presents an hourly operating workflow that considers several plausible price and demand paths rather than relying on a single expected future. It helps the station protect a small service reserve while using the remaining batteries for economical energy transactions. The reserve level, battery-aging cost, charging limits, and grid-discharge limits remain visible to and adjustable by the operator. The study is simulation based and uses one market data set and one station configuration; deployment will require local calibration and testing with actual batteries, chargers, communication delays, and market rules. The same forecast-then-decide workflow may also be useful for stationary storage, charging facilities, microgrids, and other energy assets that must act before uncertain prices or demands are revealed.
\end{abstract}}


\section{Introduction}

\IEEEPARstart{B}{attery-swapping} stations (BSSs) offer an alternative to waiting for an electric vehicle (EV) to recharge. A station stores and charges a shared inventory of standardized battery packs; when an EV arrives, its depleted pack is removed and replaced by a service-ready pack. The exchange can restore driving range rapidly, reduce charging-related waiting, and ease range anxiety, while centralized charging allows the battery inventory to be monitored and managed as a fleet rather than as isolated vehicle assets \cite{mak2013bss,zhan2022review,cui2023operation}.

That inventory gives a BSS a second role. Because many packs are connected to station chargers while they wait for customers, the station can act as grid-connected energy storage: it can charge when electricity is inexpensive or abundant, postpone charging during stressed periods, and export stored energy when market conditions justify it. This flexibility can shift demand away from critical hours and support the integration of variable energy resources. The customer-service and grid roles are nevertheless coupled. A pack discharged for an attractive market opportunity may be unavailable when the next driver arrives, whereas keeping every pack full protects service but wastes storage flexibility and may increase energy cost \cite{sarker2015bss,tan2019bss,mahoor2019least,nayak2024bss,chen2023valuation}.

Real station operation therefore follows a forecast-then-decide loop. Before scheduling the batteries, the operator must anticipate the next hours of swapping demand and electricity prices; only then can it decide when to charge, hold, or discharge each pack. Perfect future information is unavailable, so these decisions must be revised as new prices, arrivals, and battery states are observed. The forecast-then-decide challenge therefore extends well beyond BSSs, and the proposed framework can be adapted to the dynamic operation of electrolyzers, stationary batteries, and other energy systems \cite{kumar2018stationary,kumar2019benchmarking,hoang2024probabilistic}.

The practical challenge is that a single expected trajectory does not describe the futures against which the station must remain viable. An average demand forecast can hide a surge that exhausts the ready-pack inventory, and an average price forecast can miss the spike that changes whether charging or grid discharge is economical. Conversely, poorly focused uncertainty scenarios can make the controller protect against too many unlikely futures and leave the battery fleet underused. A useful forecasting layer must therefore identify both the expected temporal pattern and the deviations that matter to the station's cost and service decisions \cite{gneiting2007strict,gneiting2014probabilistic,hong2016gefcom,nowotarski2018probepf}.

To address this operating problem, we develop a forecast-aware, scenario-based stochastic model predictive control (SMPC) framework. A lightweight DLinear forecaster estimates multi-hour price and demand trajectories, and an uncertainty-quantification (UQ) layer converts them into a finite set of temporally structured scenarios \cite{lim2021tft,oreshkin2020nbeats,zeng2023dlinear,challu2023nhits,nie2023patchtst,liu2016svgd}. The controller evaluates charging, grid discharge, swapping service, terminal readiness, and battery degradation across those possible futures, implements the first decision, and repeats the process when the next observation arrives. Forecasting and station control remain modular, but their interface is designed around the information needed for operation rather than prediction accuracy alone \cite{donti2017task,wilder2019melding,elmachtoub2022smart,calafiore2005uncertain,campi2008exact}.

\begin{figure}[t]
    \centering
    \includegraphics[width=\columnwidth]{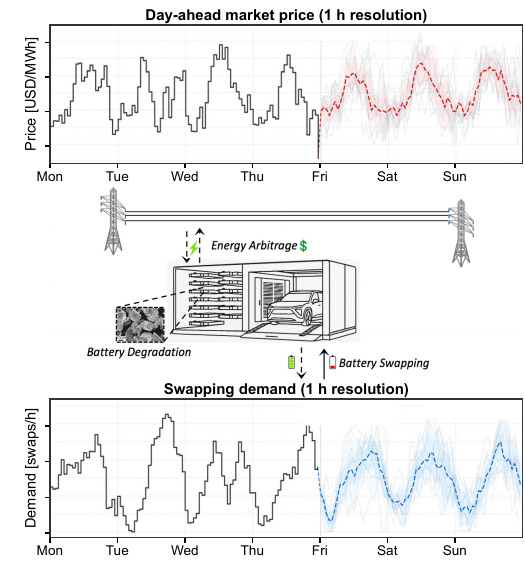}
    \caption{\revred{Operating and control framework of a battery-swapping station (BSS). The station has two coupled functions: it provides mobility service by replacing depleted EV batteries with charged batteries, and it uses the station battery inventory for energy arbitrage through grid charging and discharging. Historical price and swapping-demand data are used to forecast future trajectories, quantify their uncertainty, and generate a set of possible scenarios. A scenario-based stochastic MPC then uses these scenarios to jointly optimize battery charging, grid discharging, and swapping decisions while maintaining service readiness.}}
    \label{fig:framework}
\end{figure}

The station model reflects the physical coupling that makes the application distinctive. Swap demand removes fully charged packs from service inventory, returned depleted packs re-enter the charging pool, grid discharge competes with near-term readiness, and repeated cycling creates long-term asset cost \cite{vetter2005aging,safari2011aging,wang2011cyclelife,schmalstieg2014holistic}. We therefore include service shortfall, terminal readiness, an explicit front buffer, and battery-state propagation through a reduced-order electrochemical single-particle model (SPM). Fig.~\ref{fig:framework} summarizes how observed station data are converted into scenarios and then into hourly operating decisions.

The contributions of this work are twofold:
\begin{itemize}
    \item \textbf{BSS automation contribution:} We develop an implementable rolling-horizon controller that coordinates swapping service and grid-energy transactions using observed history and current battery states rather than perfect future information. The formulation makes service reserves, charging and discharge limits, and battery degradation explicit operating quantities.

    \item \textbf{Forecast-to-decision contribution:} We construct price--demand scenarios from lightweight DLinear forecasting and SVGD-based UQ for use in station scheduling. Matched comparisons with autoregressive forecasting, deterministic MPC, and Gaussian scenarios distinguish how mean-forecast accuracy and uncertainty representation affect BSS cost, service shortfall, and battery utilization.
\end{itemize}

The remainder of this paper is organized as follows. Section~\ref{sec:related} reviews related work. Section~\ref{sec:problem} formulates the BSS rolling energy-scheduling problem. Section~\ref{sec:forecast} presents DLinear-UQ scenario generation. Section~\ref{sec:smpc} gives the scenario-based SMPC controller. Section~\ref{sec:experiments} reports the closed-loop results, and Section~\ref{sec:conclusion} concludes the paper.

\section{Related Work}\label{sec:related}

\subsection{BSS Scheduling and EV Service Automation}
BSS planning and operation have been studied from infrastructure, power-system, and service-operations perspectives. Mak \emph{et al.} examined infrastructure planning for EVs with battery swapping, emphasizing the coupling between station location, service access, and battery inventory \cite{mak2013bss,zhan2022review,cui2023operation}. Sarker \emph{et al.} formulated operation and service scheduling for a BSS that participates in energy markets under demand and price uncertainty \cite{sarker2015bss,mahoor2019least}. Tan \emph{et al.} modeled the charging schedule of a battery charging station with dynamic inventory of fully charged batteries and quadratic degradation cost \cite{tan2019bss,chen2023valuation}. More recent work has extended BSS scheduling to stochastic demand, heterogeneous batteries, charger types, degradation, and time-of-use prices \cite{nayak2024bss,zhan2022review,cui2023operation}. In parallel, EV charging stations have become a relevant application domain for automation and predictive scheduling, including demand-charge-aware MPC \cite{yang2024tase}.

The remaining application gap is real-time operation without known future prices or customer arrivals. Many BSS scheduling models specify future profiles externally, leaving unclear how an operator should construct each rolling decision. We address this gap by treating forecasting, service inventory, and grid-energy transactions as one workflow evaluated by both cost and ready-battery availability.

\subsection{Probabilistic Energy Forecasting}
Energy forecasting has moved from point prediction toward probabilistic forecasting because power-system decisions often face asymmetric cost and reliability risks \cite{weron2014price,nowotarski2018probepf,hong2020energy}. The GEFCom2014 competition helped establish probabilistic load, price, wind, and solar forecasting as operationally important benchmarks \cite{hong2016gefcom}. Proper scoring rules and calibration--sharpness principles provide the statistical basis for evaluating distributional forecasts \cite{gneiting2007strict,gneiting2014probabilistic}. Modern probabilistic methods include quantile regression \cite{koenker1978regression,taieb2016smartmeter}, deep ensembles \cite{lakshminarayanan2017deep}, autoregressive neural density models such as DeepAR \cite{salinas2020deepar}, multi-horizon neural forecasting architectures \cite{lim2021tft,oreshkin2020nbeats,challu2023nhits,nie2023patchtst}, and particle-based Bayesian approximations such as Stein variational gradient descent (SVGD) \cite{liu2016svgd}.

{DLinear and related linear decomposition models have also shown that simple trend--seasonal structures can be strong long-horizon forecasting baselines \cite{oreshkin2020nbeats,zeng2023dlinear,challu2023nhits,nie2023patchtst}. We use DLinear because its direct multi-horizon structure is computationally light enough for repeated station operation and transparent enough to inspect when forecasts change. Its outputs are converted into scenario matrices so that forecast uncertainty can influence ready-pack inventory and energy transactions rather than remaining a separate statistical report.}

\subsection{Control-Aware Prediction and Stochastic MPC}
The forecast-then-optimize paradigm has motivated learning and evaluation criteria that reflect downstream decision quality rather than prediction error alone. Smart predict-then-optimize training links prediction models to optimization loss \cite{elmachtoub2022smart,wilder2019melding}, and task-based learning has demonstrated the value of training stochastic models through downstream energy objectives \cite{donti2017task}. Our approach follows the same control-aware principle but keeps forecasting and control modular: the forecaster is not trained end-to-end through the SMPC solver, but its output is shaped into a controller-facing scenario representation.

MPC repeatedly solves a finite-horizon problem and implements the first action, making it well suited to rolling energy management \cite{kouvaritakis2016mpc,rawlings2017mpc}. SMPC extends this idea to uncertain systems by optimizing expected or risk-aware costs and enforcing probabilistic or scenario-based constraints \cite{calafiore2005uncertain,campi2008exact,mesbah2016smpc,schildbach2014scenario}. Kumar \emph{et al.} benchmarked deterministic, stochastic, and perfect-information MPC for stationary batteries and showed that mean forecasts can cause constraint violations, while fixed reserve back-offs trade asset utilization for robustness \cite{kumar2018stationary,kumar2019benchmarking}. Recent power-system work has integrated probabilistic forecasting with stochastic nonlinear MPC for systems with renewable uncertainty and storage \cite{hewing2020learning,hoang2024probabilistic}. {We adapt these ideas to the coupled service-and-storage operation of a BSS, where uncertainty affects both energy cost and the immediate ability to serve a driver.}

\section{BSS Rolling Scheduling Problem}\label{sec:problem}

\subsection{System and Information State}
The BSS operates over hourly periods $t=0,1,\ldots,T_{\mathrm{sim}}-1$ with sampling interval $\Delta t$. At time $t$, the exogenous vector is
\begin{equation}
    \bm{\xi}_t=[d_t,\lambda_t]^\top,
\end{equation}
where $d_t$ is swapping demand and $\lambda_t$ is the electricity price. The controller observes the information set
\begin{equation}
    \calI_t=\{\bm{\xi}_{t-L},\ldots,\bm{\xi}_{t-1},\bm{X}_t\},
\end{equation}
where $L$ is the forecast look-back length and $\bm{X}_t=\{\bm{x}_{i,t}\}_{i=1}^{N_b}$ is the set of battery states for the $N_b$ station packs. Each state is
\begin{equation}
    \bm{x}_{i,t}=[\SOC_{i,t},\bm{z}_{i,t}^{\top}]^\top,
\end{equation}
where $\bm{z}_{i,t}$ collects reduced electrochemical and aging states used by the SPM rollout.

\begin{table}[t]
\centering
\caption{Notation used in the BSS forecast-to-control formulation.}
\label{tab:notation}
\footnotesize
\setlength{\tabcolsep}{3.0pt}
\renewcommand{\arraystretch}{1.04}
\begin{tabular}{@{}>{\raggedright\arraybackslash}p{0.24\columnwidth}>{\raggedright\arraybackslash}p{0.70\columnwidth}@{}}
\toprule
Symbol & Meaning \\
\midrule
$t,h$ & current time and look-ahead indices \\
$L,H$ & forecast look-back length and control horizon \\
$S_t,N_b$ & number of scenarios and station battery packs \\
$d_t,\lambda_t$ & realized swapping demand and electricity price \\
$\bm{X}_t$ & set of station battery states \\
$\SOC_{\mathrm{sw}}$ & SOC threshold for a service-ready pack \\
$\rho,N_f(\rho)$ & front-buffer ratio and protected-pack count \\
$p^{\mathrm{ch}}_{i,t+h|t}$ & planned charging power of pack $i$ \\
$p^{\mathrm{dis}}_{i,t+h|t}$ & planned grid-discharging power of pack $i$ \\
$u_{i,t+h|t},\ \bm{U}_t$ & signed pack power and control plan \\
$X_t^s$ & future demand--price trajectory in scenario $s$ \\
$\pi_t^s$ & probability of scenario $s$ \\
$\bm{D}_t,\bm{\Lambda}_t$ & demand and price scenario matrices \\
$a_{t+h|t}^{s}$ & service-ready inventory in scenario $s$ \\
$n_{\mathrm{sf},t+h|t}^{s}$ & service shortfall in scenario $s$ \\
$f_{\mathrm{SPM}}$ & single-particle-model state transition \\
$\Phi_{\mathrm{svc}}$ & FIFO swapping-service and return map \\
$c_{\mathrm{sf}},c_{\mathrm{deg}}$ & shortfall and degradation cost coefficients \\
\bottomrule
\end{tabular}
\end{table}

\subsection{Forecast Scenarios}
At each rolling step, the controller optimizes over a horizon $h=0,\ldots,H-1$. Scenario $s$ is a future demand--price trajectory
\begin{equation}
    X_{t}^{s}
    =\{(d_{t+h|t}^{s},\lambda_{t+h|t}^{s})\}_{h=0}^{H-1},
    \qquad s=1,\ldots,S_t,
\end{equation}
with probability weight $\pi_t^s$ and $\sum_{s=1}^{S_t}\pi_t^s=1$. Equivalently, the SMPC optimizer receives the matrices
\begin{equation}
    \bm{D}_t=[d_{t+h|t}^{s}]_{h,s}\in\R^{H\times S_t},\quad
    \bm{\Lambda}_t=[\lambda_{t+h|t}^{s}]_{h,s}\in\R^{H\times S_t}.
    \label{eq:scenario_matrices}
\end{equation}
These matrices are the mathematical interface between forecasting and control.

\subsection{Service Buffer and Battery-State Dynamics}
A battery pack is considered service-ready when its state of charge satisfies
$\mathrm{SOC}_{i,t}\ge \mathrm{SOC}_{\mathrm{sw}}$. The station uses a strict front buffer to reserve service-critical inventory while leaving the remaining packs available for energy arbitrage:
\begin{equation}
    N_f(\rho)=\lfloor \rho N_b \rfloor,\qquad N_a(\rho)=N_b-N_f(\rho),
    \label{eq:front_buffer}
\end{equation}
where $\rho\in[0,1]$ is the buffer ratio and $N_a(\rho)$ is the remaining energy-arbitrage-capable pool. The protected set $\calB_f(t)$ contains the $N_f(\rho)$ ready packs reserved for near-term swapping service. In the strict-buffer implementation,
\begin{equation}
    p_{i,t+h|t}^{\mathrm{dis}}=0,\qquad i\in\calB_f(t+h|t),
    \label{eq:front_buffer_discharge}
\end{equation}
so the controller cannot earn energy-sale revenue by depleting the protected service inventory.

For each battery, the nonnegative charging and discharging controls are
\begin{equation}
    p_{i,t+h|t}^{\mathrm{ch}}\ge 0,\qquad
    p_{i,t+h|t}^{\mathrm{dis}}\ge 0,
\end{equation}
with signed net power
\begin{equation}
    u_{i,t+h|t}=p_{i,t+h|t}^{\mathrm{ch}}-p_{i,t+h|t}^{\mathrm{dis}}.
\end{equation}
Here, SPM denotes a reduced-order electrochemical \emph{single-particle model}; related degradation-aware BSS control has likewise used SPM-based electrochemical aging as the physical basis for MPC scheduling \cite{li2026degradation}. It is not an optimization or evolutionary algorithm. The phrase ``SPM-consistent'' means only that the battery-state and degradation updates are propagated by this physical model. The pre-service transition is
\begin{equation}
    \widetilde{\bm{x}}_{i,t+h+1|t}^{s}
    =f_{\mathrm{SPM}}\!\left(\bm{x}_{i,t+h|t}^{s},
    u_{i,t+h|t},\Delta t\right).
    \label{eq:spm_transition}
\end{equation}
The ready inventory before service is
\begin{equation}
    a_{t+h+1|t}^{s}
    =\sum_{i=1}^{N_b}\mathbf{1}\{\widetilde{\SOC}_{i,t+h+1|t}^{s}
    \ge \SOC_{\mathrm{sw}}\},
\end{equation}
and the service shortfall is
\begin{equation}
    n_{\mathrm{sf},t+h|t}^{s}
    =\left[d_{t+h|t}^{s}-a_{t+h|t}^{s}\right]_+.
\end{equation}
After swaps are served, depleted packs are returned and the post-service station state is compactly written as
\begin{equation}
    \bm{X}_{t+h+1|t}^{s}
    =\Phi_{\mathrm{svc}}\!\left(\widetilde{\bm{X}}_{t+h+1|t}^{s},
    d_{t+h|t}^{s}\right).
\end{equation}

\section{DLinear-UQ Scenario Generation}\label{sec:forecast}

\subsection{Forecast Backbone}
Let $\bm{Y}_{t-L:t-1}\in\R^{L\times 2}$ be the historical demand--price window. DLinear decomposes the input into moving-average trend and residual components,
\begin{align}
    \bm{Y}_{t-L:t-1}^{\mathrm{tr}}
    &=\mathrm{MA}_k(\bm{Y}_{t-L:t-1}),\\
    \bm{Y}_{t-L:t-1}^{\mathrm{se}}
    &=\bm{Y}_{t-L:t-1}-\bm{Y}_{t-L:t-1}^{\mathrm{tr}},
\end{align}
and predicts the multi-horizon mean trajectory as
\begin{equation}
    \widehat{\bm{M}}_t
    =W_{\mathrm{tr}}\bm{Y}_{t-L:t-1}^{\mathrm{tr}}
     +W_{\mathrm{se}}\bm{Y}_{t-L:t-1}^{\mathrm{se}}+\bm{B},
    \label{eq:dlinear}
\end{equation}
where $\widehat{\bm{M}}_t\in\R^{H\times 2}$. Equation~\eqref{eq:dlinear} is only the forecasting backbone; the controller receives scenarios around this forecast rather than the mean trajectory alone.

\subsection{Uncertainty Quantification}
Two UQ routes can be used. In the particle route, $M$ DLinear particles are maintained and diversified by an SVGD update,
\begin{equation}
    \theta_i \leftarrow \theta_i+\epsilon\frac{1}{M}
    \sum_{j=1}^{M}\left[
    k(\theta_j,\theta_i)\nabla_{\theta_j}\log p(\theta_j|\mathcal{D})
    +\nabla_{\theta_j}k(\theta_j,\theta_i)\right].
\end{equation}
The particle predictive distribution is represented by the empirical measure over particle forecasts. In the quantile route, conditional quantiles are estimated by pinball loss and calibrated so that samples preserve temporal dependence across the horizon. The resulting finite empirical distribution is exported as $\calS_t=\{(X_t^s,\pi_t^s)\}_{s=1}^{S_t}$ and the matrices in \eqref{eq:scenario_matrices}.

The particle and quantile routes are not intended to replace standard probabilistic-forecast validation. They are used here to construct finite future trajectories with enough temporal coherence to serve as a meaningful control input for a rolling BSS scheduler.

\section{Scenario-Based SMPC}\label{sec:smpc}

\subsection{Optimization Problem}
Given $\calS_t$ and $\bm{X}_t$, the SMPC problem is
\begin{subequations}
\begin{align}
    \bm{U}_t^\star
    &\in \argmin_{\bm{U}_t}
    \sum_{s=1}^{S_t}\pi_t^s J_t^s(\bm{U}_t) \label{eq:smpc}\\
    \mathrm{s.t.}\quad
    &\bm{x}_{i,t+h+1|t}^{s}=f_{\mathrm{SPM}}
    (\bm{x}_{i,t+h|t}^{s},u_{i,t+h|t},\Delta t),\\
    &\bm{X}_{t+h+1|t}^{s}
    =\Phi_{\mathrm{svc}}(\widetilde{\bm{X}}_{t+h+1|t}^{s},
    d_{t+h|t}^{s}),\\
    &\SOC_{\min}\le \SOC_{i,t+h|t}^{s}\le \SOC_{\max},\\
    &0\le p_{i,t+h|t}^{\mathrm{ch}}\le \bar p_i^{\mathrm{ch}},
      \quad
      0\le p_{i,t+h|t}^{\mathrm{dis}}\le \bar p_i^{\mathrm{dis}},\\
    &\sum_i p_{i,t+h|t}^{\mathrm{ch}}\le \bar P^{\mathrm{ch}},
      \quad
      \sum_i p_{i,t+h|t}^{\mathrm{dis}}\le \bar P^{\mathrm{dis}},\\
    &p_{i,t+h|t}^{\mathrm{dis}}=0,\quad i\in\calB_f(t+h|t),
\end{align}
\end{subequations}
for all admissible batteries, horizons, and scenarios. The control sequence is non-anticipative: the same first action is evaluated against every scenario, and only that action is applied before the next observation arrives.

For scenario $s$, the stage cost is
\begin{align}
    \ell_{t+h|t}^{s}
    &=\lambda_{t+h|t}^{s}(P_{t+h|t}^{\mathrm{ch}}
    -\eta_{\mathrm{dis}}P_{t+h|t}^{\mathrm{dis}})\Delta t \notag\\
    &\quad+c_{\mathrm{sf}}n_{\mathrm{sf},t+h|t}^{s}
    +c_{\mathrm{deg}}\sum_i \Delta q_{i,t+h|t}^{s},
    \label{eq:stage_cost}
\end{align}
where $P^{\mathrm{ch}}$ and $P^{\mathrm{dis}}$ are aggregate station charging and grid-discharging power. The terminal value rewards end-of-horizon readiness,
\begin{equation}
    V_f(\bm{X}_{t+H|t}^{s})=
    -w_{\mathrm{soc}}\sum_i \SOC_{i,t+H|t}^{s}
    -w_a a_{t+H|t}^{s}.
\end{equation}
The scenario horizon cost is
\begin{equation}
    J_t^s(\bm{U}_t)=\sum_{h=0}^{H-1}\ell_{t+h|t}^{s}
    +V_f(\bm{X}_{t+H|t}^{s}).
\end{equation}

\begin{algorithm}[t]
\caption{DLinear-UQ-Driven Scenario SMPC for BSS Operation}
\label{alg:smpc}
\begin{algorithmic}[1]
\Require Initial battery states $\bm{X}_0$, historical observations, horizon $H$, scenario count $S$, buffer ratio $\rho$
\For{$t=0,1,\ldots,T_{\mathrm{sim}}-1$}
    \State Build $\calI_t$ from recent price--demand observations and current battery states.
    \State Generate $\calS_t$ and the matrices $\bm{D}_t$, $\bm{\Lambda}_t$ using DLinear-UQ.
    \State Form the strict front buffer $\calB_f(t)$ from ready packs and $\rho$.
    \State Solve \eqref{eq:smpc} using battery-state rollouts generated by $f_{\mathrm{SPM}}$.
    \State Apply the first charging/discharging decision.
    \State Execute charging, grid discharge, FIFO swapping service, and battery-state updates under realized $(d_t,\lambda_t)$.
\EndFor
\end{algorithmic}
\end{algorithm}

\section{Experimental Design and Results}\label{sec:experiments}

{The closed-loop study evaluates whether the proposed operating framework helps a BSS use its battery inventory economically while preserving reliable swapping service. The analysis addresses four application-facing questions:}
\begin{itemize}
    \setlength{\itemsep}{0.15em}
    \setlength{\topsep}{0.25em}
    \item \textbf{RQ1:} \textit{How well do the forecasts describe the price and swapping demand that the station will face over the next 24 hours?}
    \item \textbf{RQ2:} \textit{How do forecast accuracy and uncertainty representation affect the station's operating cost and service shortfall?}
    \item \textbf{RQ3:} \textit{How does uncertainty information change the station's charging, grid-discharge, and inventory decisions?}
    \item \textbf{RQ4:} \textit{How much battery inventory should be protected for near-term swapping service?}
\end{itemize}

\subsection{Data, Baselines, and Metrics}
\begin{figure}[t]
    \centering
    \subfloat[PJM day-ahead grid price over 120 days.]{
        \includegraphics[width=0.97\columnwidth]{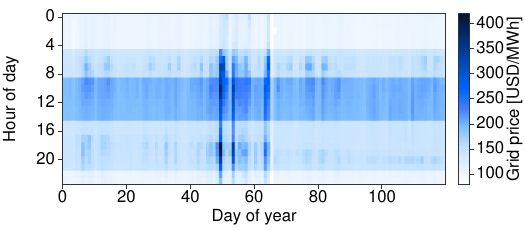}}
    \vspace{0.35em}
    \subfloat[BSS swapping demand over 120 days.]{
        \includegraphics[width=0.97\columnwidth]{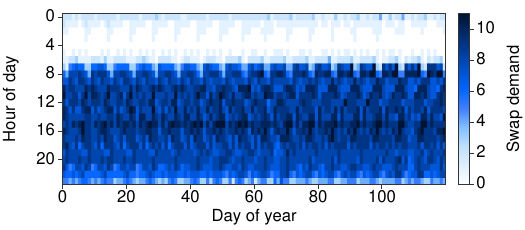}}
    \caption{{Hourly price and demand structure used in the closed-loop evaluation. The BSS controller must manage both market volatility and service-load concentration.}}
    \label{fig:data_heatmaps}
\end{figure}

{The evaluation covers 120 days at hourly resolution using PJM day-ahead prices and stochastic BSS swapping demand. Both deterministic MPC and SMPC use $H=24$~h, and DLinear uses $L=72$~h. Fig.~\ref{fig:data_heatmaps} shows intermittent price spikes and daytime/evening demand concentrations, requiring the station to anticipate demand peaks while exploiting price variation.}

{
At each decision time, the forecasting module generates $S=200$ joint price--demand scenarios. AR(24) is the autoregressive baseline and DLinear is the decomposition-based forecasting backbone. The seven reported controller configurations are:
\begin{itemize}
    \setlength{\leftmargini}{1.25em}
    \setlength{\labelsep}{0.35em}
    \setlength{\itemsep}{0pt}
    \setlength{\parskip}{0pt}
    \item {\textbf{MPC with perfect information:} uses the realized price and swapping demand for the next 24 hours in a deterministic optimization. This perfect future information yields the oracle lower bound on cost under the same model and constraints.}
    \item {\textbf{Deterministic MPC:} solves a deterministic optimization using only the predicted mean price and demand, without confidence intervals or uncertainty scenarios.}
    \item {\raggedright\textbf{Scenario-based stochastic MPC:} uses the predicted mean and UQ-generated price--demand scenarios in a  stochastic optimization. Immediate charging, discharging, and swapping decisions are shared across scenarios; future recourse is scenario dependent. The first action is implemented and the problem is solved again hourly.}
\end{itemize}
MPC with perfect information is only an oracle reference; performance claims concern the implementable deterministic and stochastic controllers.

Robustness and economic optimality are reported separately. $N_{\mathrm{sf}}$ counts hourly decisions with a positive shortfall, and $R_{\mathrm{sf}}=100(1-N_{\mathrm{sf}}/2880)\%$ is the shortfall-free-hour rate. Economic performance is
\begin{equation}
\begin{aligned}
J_{\mathrm{tr}}&=\sum_t\lambda_t(P_t^{\mathrm{ch}}-\eta_{\mathrm{dis}}P_t^{\mathrm{dis}})\Delta t,\\
J_{\mathrm{deg}}&=\sum_t c_{\mathrm{deg}}\sum_i\Delta q_{i,t},\quad
J_{\mathrm{final}}=J_{\mathrm{tr}}+J_{\mathrm{deg}},
\end{aligned}\label{eq:reported_costs}
\end{equation}
where $J_{\mathrm{tr}}$ is net grid-transaction cost after energy-sale revenue and $J_{\mathrm{deg}}$ is degradation cost. No post-hoc service cost is added; robustness is reported independently.
}

\subsection{{Price and Demand Information Available to the Station}}
\begin{figure}[t]
    \centering
    \includegraphics[width=0.96\columnwidth]{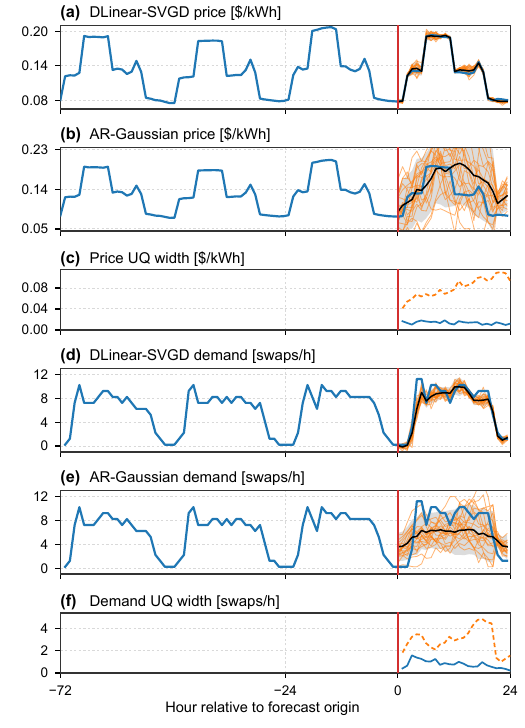}
    \caption{\revred{Mean-forecast and uncertainty comparison over 24 hours. Panels (a,b) show price and panels (d,e) demand; observations are blue, scenario means black, scenarios orange, and uncertainty bands gray. Panels (c,f) compare band widths (DLinear-SVGD: blue solid; AR-Gaussian: orange dashed). DLinear follows the observations more closely, while SVGD produces a narrower, more concentrated scenario set in this example.}}
    \label{fig:forecast_uq_integrated}
\end{figure}

\revred{RQ1 asks whether the station receives a useful picture of the next operating day. That picture has two parts: the expected price and demand trajectory, and plausible deviations around it. Fig.~\ref{fig:forecast_uq_integrated} compares both parts for price in panels (a--c) and demand in panels (d--f), using the same $L=72$~h history and $H=24$~h horizon.}

\revred{For the mean-prediction part of RQ1, DLinear reduces price RMSE from $33.6$ to $3.6$ and MAE from $28.1$ to $2.8$ over the original $T=72$~h scoring window. Demand RMSE falls from $3.2$ to $1.1$ and MAE from $2.9$ to $0.8$. The displayed $H=24$~h segment shows the same pattern: the DLinear mean follows the observations more closely.}

\revred{For the uncertainty part of RQ1, useful scenarios must capture relevant deviations without becoming so wide that every decision is conservative. Fig.~\ref{fig:forecast_uq_integrated} shows narrower SVGD bands and scenarios concentrated closer to the realized trajectory than the Gaussian baseline. In this case, SVGD therefore gives the controller a more precise representation of future uncertainty. The next subsection shows that these focused scenarios reduce both cost and service shortfall when used by stochastic MPC.}

\subsection{{Station Cost and Service Reliability}}
\revred{RQ2 asks whether better information produces better station operation. Table~\ref{tab:main_results} reports matched 120-day cost and service results, with MPC with perfect information included only as an unavailable oracle reference. DLinear-SVGD has the lowest implementable final cost, $11031.97$, and 26 shortfall hours ($99.10\%$ shortfall-free). Relative to deterministic DLinear MPC, it reduces final cost by 1.2\% and shortfall hours by 80.7\%.}

\begin{table}[t]
\centering
\caption{\revred{Closed-loop service robustness and final economic cost over 120 days. Lower $N_{\mathrm{sf}}$ and $J_{\mathrm{final}}$, and higher $R_{\mathrm{sf}}$, are better. MPC with perfect information uses the realized next-24-hour trajectories and is reported only as an oracle lower bound on cost.}}
\label{tab:main_results}
\scriptsize
\setlength{\tabcolsep}{3.0pt}
\resizebox{\columnwidth}{!}{%
{\begin{tabular}{lrrr}
\toprule
Controller & $N_{\mathrm{sf}}$ & $R_{\mathrm{sf}}$ [\%] & $J_{\mathrm{final}}$ \\
\midrule
{MPC with perfect information} & 0 & 100.00 & 9215.79 \\
\midrule
DLinear DMPC & 135 & 95.31 & 11169.12 \\
AR DMPC & 216 & 92.50 & 11523.09 \\
\midrule
DLinear-SVGD SMPC & 26 & 99.10 & \textbf{11031.97} \\
DLinear-Gaussian SMPC & 31 & 98.92 & 11113.04 \\
AR-SVGD SMPC & 1 & 99.97 & 12116.18 \\
AR-Gaussian SMPC & \textbf{0} & 100.00 & 12120.93 \\
\bottomrule
\end{tabular}}
}
\vspace{0.35em}
\begin{minipage}{\columnwidth}
\footnotesize{\textit{Note:} $R_{\mathrm{sf}}=100(1-N_{\mathrm{sf}}/2880)\%$ measures shortfall-free hours, not per-request fulfillment. $J_{\mathrm{final}}=J_{\mathrm{tr}}+J_{\mathrm{deg}}$ contains no post-hoc service penalty. Fig.~\ref{fig:cumulative_cost_main} reports the component trajectories for four representative controllers. Costs use common simulation units.}
\end{minipage}
\end{table}

\revred{The matched rows isolate both RQ2 effects. With mean forecasts only, DLinear lowers cost by 3.1\% and shortfall hours by 37.5\% relative to AR. Holding DLinear fixed, SVGD improves on Gaussian scenarios by 0.7\% in cost and 16.1\% in shortfall hours. The AR stochastic controllers reach zero or one shortfall hour only by retaining excessive inventory and paying 8.9--9.0\% more than DLinear-SVGD. Thus better UQ places protection around operationally relevant futures instead of making every decision conservative.}

\subsection{{Station Decisions Under Uncertainty}}
\begin{figure}[t]
    \centering
    \includegraphics[width=0.94\columnwidth]{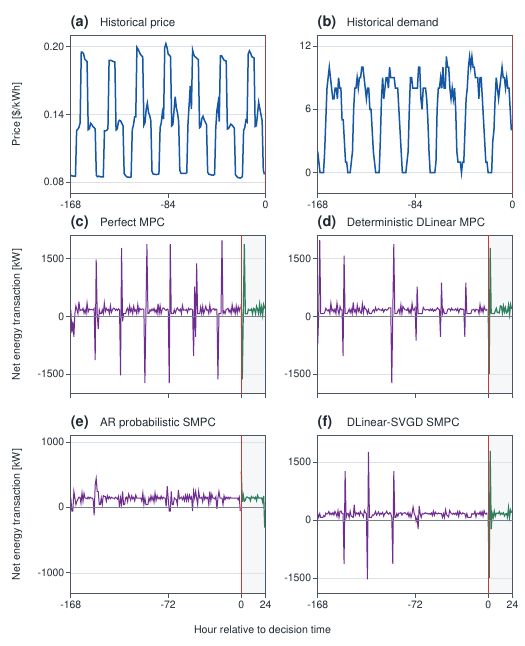}
    \caption{\revred{Representative policy snapshot under a common operating context. Panels (a,b) show realized price and demand over the preceding 168 hours. Panels (c--f) show implemented net-transaction histories and current 24-hour plans for MPC with perfect information (labeled Perfect MPC), deterministic DLinear MPC, AR probabilistic SMPC, and DLinear-SVGD SMPC. Purple denotes implemented history, green denotes the current open-loop plan, and red marks the decision time; only the first green action is implemented. Positive values are grid purchases and negative values are grid exports.}}
    \label{fig:charging_forecast}
\end{figure}

\revred{RQ3 concerns how forecast information changes what the station actually does. Fig.~\ref{fig:charging_forecast} compares behavior, while Table~\ref{tab:main_results} gives the resulting performance. Panels (a,b) show the common 168-hour price and demand history; panels (c--f) show each implemented transaction history and current 24-hour plan. Only the first planned action is implemented. The practical question is whether an energy-market opportunity remains worthwhile after accounting for the risk of leaving too few batteries ready for customers.}

\revred{MPC with perfect information is the oracle timing reference because its transactions use the realized next 24 hours. Deterministic DLinear MPC also reacts to recurring prices, but its single mean removes demand and price tails. It is therefore \emph{tail-blind}: a profitable mean-based transaction can leave too little ready inventory during a demand surge, consistent with its 135 shortfall hours.}

\revred{AR probabilistic SMPC completes the RQ3 contrast: its diffuse scenarios keep transactions near zero, yielding near-zero shortfall at the highest implementable cost. DLinear-SVGD instead tests DLinear's price opportunities against demand and price tails. This targeted protection reduces shortfall hours from 135 to 26 while also lowering cost, without imposing uniform conservatism.}

\subsection{{Service Buffer and Battery Utilization}}
\begin{figure}[t]
    \centering
    \includegraphics[width=\columnwidth]{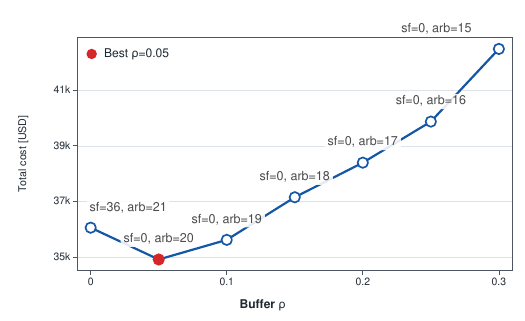}
    \caption{{Strict front-buffer sensitivity for DLinear-SVGD SMPC with a fixed 21-battery pool. Enlarged labels report shortfall hours (sf) and arbitrage-capable packs (arb) at each $\rho$. A small buffer removes shortfalls without over-constraining the arbitrage pool; excessive buffer ratios increase total cost.}}
    \label{fig:rho_sweep}
\end{figure}

\revred{RQ4 asks how much inventory should be protected for customers rather than exposed to energy arbitrage. Fig.~\ref{fig:rho_sweep} makes this operating trade-off explicit for 21 batteries. With $\rho=0$, all packs can arbitrage but 36 shortfalls occur. A small buffer, $\rho=0.05$, removes shortfalls and minimizes cost; larger buffers remain reliable but reduce the arbitrage pool and increase cost.}

\begin{figure}[t]
    \centering
    \includegraphics[width=\columnwidth]{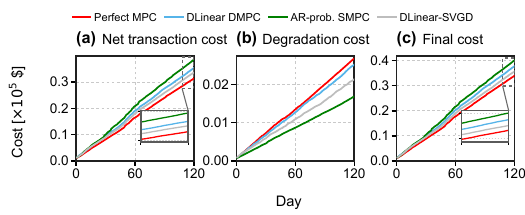}
    \caption{{Horizontal cumulative 120-day cost comparison at $\rho=0.05$: (a) net grid-transaction cost after energy-sale revenue, $J_{\mathrm{tr}}$; (b) degradation cost, $J_{\mathrm{deg}}$; and (c) final cost, $J_{\mathrm{final}}=J_{\mathrm{tr}}+J_{\mathrm{deg}}$. Dashed boxes in (a) and (c) mark days 108--120, and straight connectors lead to the corresponding terminal insets while preserving the complete trajectories.}}
    \label{fig:cumulative_cost_main}
\end{figure}

\revred{RQ4 thus supports a small, interpretable buffer rather than a large reserve. Fig.~\ref{fig:cumulative_cost_main} decomposes the 120-day results at $\rho=0.05$. AR probabilistic SMPC has the lowest degradation but the largest transaction cost, indicating battery under-utilization. DLinear-SVGD accepts moderate cycling when justified by market value and closes part of the oracle cost gap.}

\subsection{Discussion}\label{sec:discussion}

\subsubsection{{Operational Value of Anticipating Price and Demand}}
\revred{The results show how information changes the use of the station's battery inventory. A more accurate mean forecast helps the operator locate hours in which charging or grid discharge is economically attractive. The scenario set then tests those opportunities against possible demand surges and price deviations. This prevents the deterministic controller from committing too much ready inventory under one expected future, while avoiding the blanket conservatism produced by diffuse scenarios. For a BSS, the benefit is concrete: uncertainty is used to protect customer service when needed, not to keep every battery idle at all times.}

\revred{The improvements have two operational meanings. Fewer shortfall hours indicate that more drivers encounter a ready pack, while lower transaction cost shows that reliability is not obtained by charging every battery early and holding it idle. Grid discharge can monetize unused capacity and selective charging can avoid expensive hours, but both are valuable only when enough inventory remains for mobility service. Cost and shortfall must therefore be evaluated together.}

\subsubsection{{Deployment and Broader Energy Applications}}
\revred{The controller requires measured price and demand history together with current battery states; it does not rely on the perfect future used by the oracle benchmark. Operators retain direct control over the protected inventory, charging and grid-discharge limits, service penalties, and degradation cost. Before deployment, these quantities and the forecast scenarios must be calibrated with local station data, charger capabilities, communication delays, and market rules. The study is simulation based, and hardware-in-the-loop and field trials are needed to test model mismatch and real-time reliability. The same forecast-then-decide architecture can support stationary storage, microgrids, flexible charging, and electricity-driven production whenever an energy asset must be scheduled before uncertain prices or demands are known.}

\section{Conclusion}\label{sec:conclusion}
Battery swapping can provide rapid EV service and flexible grid-connected storage, but these two roles compete for the same battery inventory under uncertain demand and prices. This paper developed a forecast-aware rolling controller that coordinates service readiness, charging, grid discharge, and degradation without assuming perfect future information. The matched experiments show that improved mean forecasting makes energy transactions more selective, while a focused uncertainty representation improves the cost--service balance relative to Gaussian scenarios and deterministic mean MPC. The resulting framework makes the service buffer and physical operating limits auditable by station operators. Future work will add field-calibrated uncertainty evaluation, multi-station coordination, online recalibration, distributional robustness, and hardware-in-the-loop validation.


\begin{thebibliography}{99}
\fontsize{8.0}{8.6}\selectfont

\bibitem{mak2013bss}
H.-Y. Mak, Y. Rong, and Z.-J. M. Shen, ``Infrastructure planning for electric vehicles with battery swapping,'' \emph{Management Science}, vol. 59, no. 7, pp. 1557--1575, 2013.

\bibitem{zhan2022review}
W. Zhan, Z. Wang, L. Zhang, P. Liu, D. Cui, and D. G. Dorrell, ``A review of siting, sizing, optimal scheduling, and cost-benefit analysis for battery swapping stations,'' \emph{Energy}, vol. 258, Art. no. 124723, 2022.

\bibitem{cui2023operation}
D. Cui, Z. Wang, P. Liu, S. Wang, D. G. Dorrell, X. Li, and W. Zhan, ``Operation optimization approaches of electric vehicle battery swapping and charging station: A literature review,'' \emph{Energy}, vol. 263, Art. no. 126095, 2023.

\bibitem{sarker2015bss}
M. R. Sarker, H. Pandzic, and M. A. Ortega-Vazquez, ``Optimal operation and services scheduling for an electric vehicle battery swapping station,'' \emph{IEEE Transactions on Power Systems}, vol. 30, no. 2, pp. 901--910, Mar. 2015.

\bibitem{tan2019bss}
X. Tan, G. Qu, B. Sun, N. Li, and D. H. K. Tsang, ``Optimal scheduling of battery charging station serving electric vehicles based on battery swapping,'' \emph{IEEE Transactions on Smart Grid}, vol. 10, no. 2, pp. 1372--1384, Mar. 2019.

\bibitem{mahoor2019least}
M. Mahoor, Z. S. Hosseini, and A. Khodaei, ``Least-cost operation of a battery swapping station with random customer requests,'' \emph{Energy}, vol. 172, pp. 913--921, 2019.

\bibitem{nayak2024bss}
D. S. Nayak and S. Misra, ``An operational scheduling framework for electric vehicle battery swapping station under demand uncertainty,'' \emph{Energy}, vol. 286, Art. no. 130219, 2024.

\bibitem{chen2023valuation}
X. Chen, Y. Yang, J. Wang, J. Song, and G. He, ``Battery valuation and management for battery swapping station,'' \emph{Energy}, vol. 279, Art. no. 128120, 2023.

\bibitem{kumar2018stationary}
R. Kumar, M. J. Wenzel, M. J. Ellis, M. N. ElBsat, K. H. Drees, and V. M. Zavala, ``A stochastic model predictive control framework for stationary battery systems,'' \emph{IEEE Transactions on Power Systems}, vol. 33, no. 4, pp. 4397--4406, July 2018.

\bibitem{kumar2019benchmarking}
R. Kumar, J. Jalving, M. J. Wenzel, M. J. Ellis, M. N. ElBsat, K. H. Drees, and V. M. Zavala, ``Benchmarking stochastic and deterministic MPC: A case study in stationary battery systems,'' \emph{AIChE Journal}, vol. 65, no. 7, Art. no. e16551, 2019.

\bibitem{hoang2024probabilistic}
K. T. Hoang, C. A. Thilker, B. R. Knudsen, and L. S. Imsland, ``Probabilistic forecasting-based stochastic nonlinear model predictive control for power systems with intermittent renewables and energy storage,'' \emph{IEEE Transactions on Power Systems}, vol. 39, no. 4, pp. 5522--5534, 2024.

\bibitem{gneiting2007strict}
T. Gneiting and A. E. Raftery, ``Strictly proper scoring rules, prediction, and estimation,'' \emph{Journal of the American Statistical Association}, vol. 102, no. 477, pp. 359--378, 2007.

\bibitem{gneiting2014probabilistic}
T. Gneiting and M. Katzfuss, ``Probabilistic forecasting,'' \emph{Annual Review of Statistics and Its Application}, vol. 1, pp. 125--151, 2014.

\bibitem{hong2016gefcom}
T. Hong, P. Pinson, S. Fan, H. Zareipour, A. Troccoli, and R. J. Hyndman, ``Probabilistic energy forecasting: Global Energy Forecasting Competition 2014 and beyond,'' \emph{International Journal of Forecasting}, vol. 32, no. 3, pp. 896--913, 2016.

\bibitem{nowotarski2018probepf}
J. Nowotarski and R. Weron, ``Recent advances in electricity price forecasting: A review of probabilistic forecasting,'' \emph{Renewable and Sustainable Energy Reviews}, vol. 81, pp. 1548--1568, 2018.

\bibitem{lim2021tft}
B. Lim, S. O. Arik, N. Loeff, and T. Pfister, ``Temporal fusion transformers for interpretable multi-horizon time series forecasting,'' \emph{International Journal of Forecasting}, vol. 37, no. 4, pp. 1748--1764, 2021.

\bibitem{oreshkin2020nbeats}
B. N. Oreshkin, D. Carpov, N. Chapados, and Y. Bengio, ``N-BEATS: Neural basis expansion analysis for interpretable time series forecasting,'' in \emph{Proc. Int. Conf. Learn. Represent. (ICLR)}, 2020.

\bibitem{zeng2023dlinear}
A. Zeng, M. Chen, L. Zhang, and Q. Xu, ``Are transformers effective for time series forecasting?,'' in \emph{Proc. AAAI Conf. Artif. Intell.}, vol. 37, no. 9, pp. 11121--11128, 2023.

\bibitem{challu2023nhits}
C. Challu, K. G. Olivares, B. N. Oreshkin, F. Garza, M. Mergenthaler-Canseco, and A. Dubrawski, ``N-HiTS: Neural hierarchical interpolation for time series forecasting,'' in \emph{Proc. AAAI Conf. Artif. Intell.}, vol. 37, no. 6, pp. 6989--6997, 2023.

\bibitem{nie2023patchtst}
Y. Nie, N. H. Nguyen, P. Sinthong, and J. Kalagnanam, ``A time series is worth 64 words: Long-term forecasting with transformers,'' in \emph{Proc. Int. Conf. Learn. Represent. (ICLR)}, 2023.

\bibitem{liu2016svgd}
Q. Liu and D. Wang, ``Stein variational gradient descent: A general purpose Bayesian inference algorithm,'' in \emph{Proc. Adv. Neural Inf. Process. Syst.}, pp. 2378--2386, 2016.

\bibitem{donti2017task}
P. L. Donti, B. Amos, and J. Z. Kolter, ``Task-based end-to-end model learning in stochastic optimization,'' in \emph{Proc. Adv. Neural Inf. Process. Syst.}, pp. 5484--5494, 2017.

\bibitem{wilder2019melding}
B. Wilder, B. Dilkina, and M. Tambe, ``Melding the data-decisions pipeline: Decision-focused learning for combinatorial optimization,'' in \emph{Proc. AAAI Conf. Artif. Intell.}, vol. 33, no. 1, pp. 1658--1665, 2019.

\bibitem{elmachtoub2022smart}
A. N. Elmachtoub and P. Grigas, ``Smart `predict, then optimize','' \emph{Management Science}, vol. 68, no. 1, pp. 9--26, 2022.

\bibitem{calafiore2005uncertain}
G. C. Calafiore and M. C. Campi, ``Uncertain convex programs: Randomized solutions and confidence levels,'' \emph{Mathematical Programming}, vol. 102, no. 1, pp. 25--46, 2005.

\bibitem{campi2008exact}
M. C. Campi and S. Garatti, ``The exact feasibility of randomized solutions of uncertain convex programs,'' \emph{SIAM Journal on Optimization}, vol. 19, no. 3, pp. 1211--1230, 2008.

\bibitem{vetter2005aging}
J. Vetter, P. Novak, M. R. Wagner, C. Veit, K.-C. Moller, J. O. Besenhard, M. Winter, M. Wohlfahrt-Mehrens, C. Vogler, and A. Hammouche, ``Ageing mechanisms in lithium-ion batteries,'' \emph{Journal of Power Sources}, vol. 147, no. 1--2, pp. 269--281, 2005.

\bibitem{safari2011aging}
M. Safari and C. Delacourt, ``Aging of a commercial graphite/LiFePO$_4$ cell,'' \emph{Journal of The Electrochemical Society}, vol. 158, no. 10, pp. A1123--A1135, 2011.

\bibitem{wang2011cyclelife}
J. Wang, P. Liu, J. Hicks-Garner, E. Sherman, S. Soukiazian, M. Verbrugge, H. Tataria, J. Musser, and P. Finamore, ``Cycle-life model for graphite-LiFePO$_4$ cells,'' \emph{Journal of Power Sources}, vol. 196, no. 8, pp. 3942--3948, 2011.

\bibitem{schmalstieg2014holistic}
J. Schmalstieg, S. Kabitz, M. Ecker, and D. U. Sauer, ``A holistic aging model for Li(NiMnCo)O$_2$ based 18650 lithium-ion batteries,'' \emph{Journal of Power Sources}, vol. 257, pp. 325--334, 2014.

\bibitem{yang2024tase}
L. Yang, X. Geng, X. Guan, and L. Tong, ``EV charging scheduling under demand charge: A block model predictive control approach,'' \emph{IEEE Transactions on Automation Science and Engineering}, vol. 21, no. 2, pp. 2125--2138, Apr. 2024.

\bibitem{weron2014price}
R. Weron, ``Electricity price forecasting: A review of the state-of-the-art with a look into the future,'' \emph{International Journal of Forecasting}, vol. 30, no. 4, pp. 1030--1081, 2014.

\bibitem{hong2020energy}
T. Hong, P. Pinson, Y. Wang, R. Weron, D. Yang, and H. Zareipour, ``Energy forecasting: A review and outlook,'' \emph{IEEE Open Access Journal of Power and Energy}, vol. 7, pp. 376--388, 2020.

\bibitem{koenker1978regression}
R. Koenker and G. Bassett, Jr., ``Regression quantiles,'' \emph{Econometrica}, vol. 46, no. 1, pp. 33--50, 1978.

\bibitem{taieb2016smartmeter}
S. B. Taieb, R. Huser, R. J. Hyndman, and M. G. Genton, ``Forecasting uncertainty in electricity smart meter data by boosting additive quantile regression,'' \emph{IEEE Transactions on Smart Grid}, vol. 7, no. 5, pp. 2448--2455, 2016.

\bibitem{lakshminarayanan2017deep}
B. Lakshminarayanan, A. Pritzel, and C. Blundell, ``Simple and scalable predictive uncertainty estimation using deep ensembles,'' in \emph{Proc. Adv. Neural Inf. Process. Syst.}, pp. 6402--6413, 2017.

\bibitem{salinas2020deepar}
D. Salinas, V. Flunkert, J. Gasthaus, and T. Januschowski, ``DeepAR: Probabilistic forecasting with autoregressive recurrent networks,'' \emph{International Journal of Forecasting}, vol. 36, no. 3, pp. 1181--1191, 2020.

\bibitem{kouvaritakis2016mpc}
B. Kouvaritakis and M. Cannon, \emph{Model Predictive Control: Classical, Robust and Stochastic}. Cham, Switzerland: Springer, 2016.

\bibitem{rawlings2017mpc}
J. B. Rawlings, D. Q. Mayne, and M. M. Diehl, \emph{Model Predictive Control: Theory, Computation, and Design}, 2nd ed. Madison, WI, USA: Nob Hill Publishing, 2017.

\bibitem{mesbah2016smpc}
A. Mesbah, ``Stochastic model predictive control: An overview and perspectives for future research,'' \emph{IEEE Control Systems Magazine}, vol. 36, no. 6, pp. 30--44, Dec. 2016.

\bibitem{schildbach2014scenario}
G. Schildbach, L. Fagiano, C. Frei, and M. Morari, ``The scenario approach for stochastic model predictive control with bounds on closed-loop constraint violations,'' \emph{Automatica}, vol. 50, no. 12, pp. 3009--3018, 2014.

\bibitem{hewing2020learning}
L. Hewing, K. P. Wabersich, M. Menner, and M. N. Zeilinger, ``Learning-based model predictive control: Toward safe learning in control,'' \emph{Annual Review of Control, Robotics, and Autonomous Systems}, vol. 3, pp. 269--296, 2020.

\bibitem{li2026degradation}
R. Li, Z. Chen, Z. Zhang, R. Guo, Z. Sun, J. Yao, and J. Ma, ``Degradation-aware model predictive control for battery swapping stations under energy arbitrage,'' \emph{IEEE Transactions on Automation Science and Engineering}, accepted for publication, 2026, arXiv:2510.07902.

\end{thebibliography}
\end{document}